\documentclass[12pt]{amsart} 
\usepackage{amsmath,amssymb}
\usepackage{latexsym}
\usepackage{amssymb} 
\usepackage{latexsym}
\usepackage{delarray}

\numberwithin{equation}{section}
\theoremstyle{plain}
\newtheorem{thm}{Theorem}[section]

\title[Dispersion and asymptotic profiles]
{Dispersion and asymptotic profiles 
for Kirchhoff equations}

\author[Matsuyama and Ruzhansky]
{Tokio Matsuyama${}^\dagger$ and Michael Ruzhansky${}^*$} 
\address{ 
${}^\dagger$Department of Mathematics \endgraf 
Tokai University \endgraf 
Hiratsuka \endgraf 
Kanagawa 259-1292 \endgraf 
Japan \endgraf 
{} \endgraf 
\bigskip
${}^*$Department of Mathematics \endgraf
Imperial College London\endgraf
180 Queen's gate \endgraf 
London SW7 2AZ \endgraf 
United Kingdom} 
\thanks
{ \endgraf
 \endgraf
The second author was supported by the Leverhulme Research
Fellowship} 
\bigskip 
\email{tokio@keyaki.cc.u-tokai.ac.jp \endgraf 
m.ruzhansky@imperial.ac.uk} 

\begin{document} 

\maketitle 
\begin{abstract}
The aim of this article is to 
describe asymptotic profiles for the Kirchhoff equation, 
and to establish time decay properties and dispersive estimates for 
Kirchhoff equations. For this purpose, the method of asymptotic integration
is developed for the corresponding linear equations and representation
formulae for their solutions are obtained. 
These formulae are analysed further to obtain the
time decay rate of $L^p$--$L^q$
norms of propagators for the corresponding Cauchy problems.
\end{abstract}

\section{Introduction}
\setcounter{equation}{0} 
This article is devoted to several aspects of Kirchhoff equations or 
Kirchhoff systems, which were discussed in 
\cite{Matsuyama00,Matsuyama1,Matsuyama2}. 
In particular, we will discuss the
asymptotic profiles and dispersion properties,
or time decay of $L^p$--$L^q$ norms of propagators for some
relevant classes of hyperbolic equations.
These properties are well-known for the wave equations, but several 
aspects of Kirchhoff equations still remain far from
being understood. The global well-posedness 
of Kirchhoff equations or Kirchhoff systems is known if the data is 
sufficiently small in some suitable Sobolev spaces of $L^2$ type (see 
\cite{Manfrin,Dancona1,Dancona2,Dancona3,Greenberg,Manfrin1,Manfrin2,
Yamazaki1,Yamazaki2,Yamazaki3}).
Up to now, 
if one takes any large data from these Sobolev spaces, 
the problem of the global well-posedness is still open. 

In this article the asymptotics and the global 
well-posedness are discussed for small data. 
The first topic was developed in \cite{Matsuyama2} by relating 
the problem to 
the asymptotic behaviour of the Bessel potentials (Theorem from 
\cite{Matsuyama00} is the anoucement of \cite{Matsuyama2}). 
More precisely, the 
first author proved 
that there exists a solution which is never asymptotically free. 
Here we say 
that $u=u(t,x)$ is asymptotically free if it is asymptotically
convergent  to some 
solution of the free wave equation as the time goes to $\pm\infty$. 
From the point of view of the 
scattering theory all solutions with data 
satisfying some fast decay conditions in space variables are 
asymptotically 
free (see \cite{Ghisi,Greenberg,Yamazaki1}), while the result of 
\cite{Matsuyama2} states that if the data satisfy the 
opposite condition to 
\cite{Ghisi,Greenberg,Yamazaki1}, then the scattering theory is 
not possible. This is stated more precisely 
in Theorem \ref{thm:thm2}. 
For deriving these asymptotics, we need a delicate 
analysis of an oscillatory 
integral associated with Kirchhoff equation, which 
was introduced by Greenberg 
and Hu \cite{Greenberg} in the one dimensional case (see 
also \cite{Dancona1,Dancona2,Yamazaki1}), and we will 
develop an asymptotic 
expansion of this oscillatory integral. 

For further investigations, for example, 
such as the nonlinear scattering 
theory, the second topic is very important. 
This means that there exists 
a scattering state for Kirchhoff equations or systems with nonlinear 
perturbations, which can be discussed in the standard way but is quite 
lengthy, hence we do not touch it (see e.g., \cite{mat}). 
Quite recently, the first author obtained the 
dispersive estimates for the Kirchhoff equation 
(see \cite{Matsuyama1}), which will be introduced as Theorem \ref{thm:M}. 
The essential point of the proof relies on 
the stationary phase method together with Littman's lemma. 

Now let us give the precise formulation of Kirchhoff equations
considered problems. 
In 1883 G. Kirchhoff proposed the equation 
\begin{equation}
u_{tt}
-\bigg(1+\displaystyle\int^{L}_{0} u^{2}_x \,dx 
\bigg) u_{xx}=0
\label{Kirchhoff}
\end{equation}
for $u=u(t,x)$ on $\mathbb{R}_t \times (0,L)$ (see \cite{K}), 
which describes the nonlinear vibrations of one 
dimensional elastic strings 
having the natural length $L$. For simplicity, 
all the physical constants are normalised. 
Generalising the equation \eqref{Kirchhoff} to a multi-dimensional 
version, we can consider the Cauchy problem for $u=u(t,x)$ on 
$\mathbb{R}_t \times \mathbb{R}^{n}_x$: 
\begin{equation} 
 \partial^{2}_{t} u
-\bigg(1+\displaystyle\int_{\mathbb{R}^n}\vert\nabla u \vert^2 \,dx 
\bigg) \Delta u=0, \label{Kirchhoff1}
\end{equation}
\begin{equation}
 u(0,x)=f_0(x),\quad \partial_t u(0,x)=f_1(x), \label{Kirchhoff2}
\end{equation}
where 
$\partial_t=\frac{\partial}{\partial t}$, $\nabla=
\left(\frac{\partial}{\partial x_1},\ldots,
\frac{\partial}{\partial x_n} 
\right)$ and $\Delta$ is the standard
Laplacian in $\mathbb{R}^n$ defined by 
$\Delta=\displaystyle{\sum^{n}_{j=1}} 
\frac{\partial^2}{\partial x^{2}_{j}}$. 

Higher order nonlinear equations of Kirchhoff type are also of great 
interest, and they can be viewed as dispersion relations for Kirchhoff systems.
In particular, since higher order equations are influenced by the geometric
properties of characteristics (cf. \cite{sugi94,sugi96,sugi98}), 
in such problem it is important to know
how this phenomenon is affected by nonlinearities
(this is contrary to the $L^p$--$L^p$ estimates, see
\cite{Ruzh-survey} for a survey of such results).

Thus, let us consider the
following nonlinear equation
\begin{equation}
\widetilde{L}(t,D_t,D_x,\Vert \nabla u \Vert^2_{L^2})
=D^{m}_t u+\sum_{\underset{j \le m-1}
{\vert \nu \vert+j=m}} 
b_{\nu,j}\left(\Vert \nabla u(t,\cdot) \Vert^2_{L^2}\right) 
D^{\nu}_x D^{j}_t u=0, 
\label{Kirchhoff-Equation}
\end{equation}
for $t\not=0$,
with the initial condition 
\begin{equation}
D^{k}_t u(0,x)=f_k(x), \quad 
k=0,1,\ldots,m-1, \quad x \in \mathbb{R}^n, 
\label{Kirchhoff data}
\end{equation}
where $D_t=\frac{1}{i}\frac{\partial}{\partial t}$ and $D^{\nu}_x=
\left(\frac{1}{i}\frac{\partial}{\partial x_1} \right)^{\nu_1}
\cdots \left(\frac{1}{i}\frac{\partial}{\partial x_n} 
\right)^{\nu_n}$, 
$i=\sqrt{-1}$, for $\nu=(\nu_1,\ldots,\nu_n)$. 
We will assume that the symbol of the differential operator 
$\widetilde{L}(t,
D_t,D_x,\Vert \nabla u \Vert^2_{L^2})$ has real 
and distinct roots 
$\widetilde{\varphi}_1(t,s;\xi),\ldots,\widetilde{\varphi}_m(t,s;\xi)$ 
for $\xi \ne0$ and $0\le s \le \delta$ with $\delta>0$, i.e. 
\[
\widetilde{L}(t,\tau,\xi,s)=(\tau-\widetilde{\varphi}_1(t,s;\xi)) \cdots 
(\tau-\widetilde{\varphi}_m(t,s;\xi)), 
\] 
\[
\inf_{\underset{j \ne k}
{\vert \xi \vert=1,t \in \mathbb{R},s \in[0,\delta]}} 
\vert \widetilde{\varphi}_j(t,s;\xi)-\widetilde{\varphi}_k(t,s;\xi) \vert>0.
\]
The detail analysis of the Cauchy problem 
\eqref{Kirchhoff-Equation}--\eqref{Kirchhoff data} will be done 
in \cite{MR}, and we will consider only the Cauchy problem 
\eqref{Kirchhoff1}--\eqref{Kirchhoff2} of the second order in this article.

We conclude the introduction by fixing the 
notation used in this article. 
For $s \in \mathbb{R}$ and 
$1 \le p \le \infty$, let $\dot{L}^{p}_s=
\dot{L}^{p}_s(\mathbb{R}^n)$ and 
$L^{p}_s=L^{p}_s(\mathbb{R}^n)$ be the Riesz and 
Bessel potential spaces 
with semi-norm or norm 
\[
\Vert u \Vert_{\dot{L}^{p}_s}=\Vert \mathcal{F}^{-1}
[\vert \xi \vert^{s} \widehat{u}(\xi)] \Vert_{L^p(\mathbb{R}^n)}
\equiv \Vert \vert D \vert^s u \Vert_{L^p(\mathbb{R}^n)}, 
\]
\[ 
\Vert u \Vert_{L^{p}_s}=\Vert \mathcal{F}^{-1}
[\langle \xi \rangle^{s} \widehat{u}(\xi)] \Vert_{L^p(\mathbb{R}^n)} 
\equiv \Vert \langle D \rangle^s u \Vert_{L^p(\mathbb{R}^n)}, 
\] 
respectively. Here {} $\widehat{}$ {} denotes the Fourier transform, 
$\mathcal{F}^{-1}$ is its inverse, and 
$\langle \xi \rangle=\sqrt{1+\vert \xi \vert^2}$. 
Throughout this article, we 
fix the notation as follows: 
$$
\text{$\dot{H}^{s}=\dot{L}^{2}_s$, \quad 
$H^{s}=L^{2}_s$.} 
$$ 
We also put, for $s \ge 1$, 
\[
\dot{X}^s(\mathbb{R})=C(\mathbb{R};\dot{H}^s )\cap 
C^1(\mathbb{R};\dot{H}^{s-1}) \cap 
C^2(\mathbb{R};\dot{H}^{s-2}).
\] 
Finally we shall denote by $\mathcal{S}=\mathcal{S}
(\mathbb{R}^n)$ the Schwartz space on $\mathbb{R}^n$.

\section{Results}
\setcounter{equation}{0} 
In this section we survey the results of \cite{Matsuyama1} and 
\cite{Matsuyama2} on the Cauchy 
problem \eqref{Kirchhoff1}--\eqref{Kirchhoff2}. 
In order to state these asymptotics for the solutions to Kirchhoff 
equation, we refer to a general theorem of Yamazaki (see \cite{Yamazaki1}). For this purpose, let us introduce the set 
\[
Y_k:= \left\{ \, \{ \phi,\psi \} \in \dot{H}^{3/2} \times H^{1/2} \, ; \, 
\pmb{\vert} \left\{ \phi,\psi \right\} \pmb{\vert}_{Y_k}<\infty \, \right\}, 
\quad k>1, 
\]
where 
\begin{gather*}
\pmb{\vert} \left\{ \phi,\psi \right\} \pmb{\vert}_{Y_k} 
:= \sup_{\tau \in \mathbb{R}} (1+\vert \tau \vert)^k \left\vert 
\int\limits_{\mathbb{R}^n} \mathrm{e}^{i\tau \vert \xi \vert} 
\vert \xi \vert^3 \vert \widehat{\phi}(\xi) \vert^2 \, d \xi \right\vert \\ 
+\sup_{\tau \in \mathbb{R}} (1+\vert \tau \vert)^k \left\vert 
\int\limits_{\mathbb{R}^n} \mathrm{e}^{i\tau \vert \xi \vert} 
\vert \xi \vert \vert \widehat{\psi}(\xi) \vert^2 \, d\xi \right\vert\\ 
+\sup_{\tau \in \mathbb{R}} (1+\vert \tau \vert)^k \left\vert 
\int\limits_{\mathbb{R}^n} \mathrm{e}^{i\tau \vert \xi \vert} 
\vert \xi \vert^2 {\rm Re} \left( \widehat{\phi}(\xi) \overline{\widehat{\psi} (\xi)} 
\right) \, d \xi \right\vert. 
\end{gather*}
Then we have the following: 

\newtheorem{Ythm}{Theorem A}
\renewcommand{\theYthm}{} 
\begin{Ythm}[\cite{Yamazaki1}] Let $n \ge 1$ and $s_0 \ge \frac32$. If the 
data $u_0$, $u_1$ satisfy $u_0 \in \dot{H}^{s_0}\cap H^1$, $u_1 \in H^{s_0-1}$, and 
\begin{equation}
\text{$\delta_1:=\Vert \nabla u_0 \Vert^{2}_{L^2}+\Vert u_1 \Vert^{2}_{L^2}
+\pmb{\vert} \left\{ u_0,u_1 \right\} \pmb{\vert}_{Y_k} \ll 1$ \quad for some 
$k>1$,}
\label{smallness}
\end{equation}
then the problem \eqref{Kirchhoff1}--\eqref{Kirchhoff2} 
has a unique solution $u(t,x) \in \dot{X}^{s_0}
(\mathbb{R})$ having the following property{\rm :} there exists a constant 
$c_{\pm\infty}\equiv c_{\pm\infty}(u_0,u_1)>0$ such that 
\[
\text{$1+\Vert \nabla u(t,\cdot) \Vert^{2}_{L^2}=c^{2}_{\pm\infty}
+O \left(\vert t \vert^{-k+1} \right)$ \quad as $t \to \pm\infty$.} 
\] 
Furthermore, if \eqref{smallness} holds with $k>2$, then $c_{+\infty}
=c_{-\infty}:=c_\infty$ and each 
solution $u(t,x)\in \dot{X}^{s_0}(\mathbb{R})$ is asymptotically free in 
$\dot{H}^{\sigma}\times \dot{H}^{\sigma-1}$ for all $\sigma \in [1,s_0]$ as 
$t\to \pm \infty$, i.e., there exists a solution $v_\pm=v_\pm(t,x)\in 
\dot{X}^{\sigma}(\mathbb{R})$ of the equation 
$$\left( \partial^{2}_t-c^{2}_\infty \Delta \right) v_\pm=0 \quad 
\text{on $\mathbb{R} \times \mathbb{R}^n$}
$$ 
such that 
\[
\Vert u(t,\cdot)-v_\pm(t,\cdot)\Vert_{\dot{H}^\sigma}+
\Vert \partial_t u(t,\cdot)-\partial_t v_\pm(t,\cdot)
\Vert_{\dot{H}^{\sigma-1}}\to 0 \quad (t \to \pm\infty).
\]
\end{Ythm}

The inclusions among the classes $Y_k$ are as follows{\rm :}
\[
Y_k \subset Y_\ell \quad \text{if $k>\ell>1$, \quad and \quad 
$\mathcal{S} \subset Y_k$ \quad for all $k \in (1,n+1]$.}
\]
The latter inclusion can be shown by using the asymptotic expansion of 
oscillatory integral $I(\tilde{\vartheta}(t),0)$ which was proved in 
\cite{Matsuyama2}. 
The definition of $Y_k$ is somewhat complicated. There are some 
examples of spaces contained in $Y_k$. For more details see \cite{Matsuyama2}. 

Keeping in mind Theorem A, we have $L^p$--$L^q$ estimates: 
\begin{thm}[\cite{Matsuyama1}] \label{thm:M}
Let $n \ge 2$ and let 
$1<p \le 2 \le q<+\infty$ and $\frac{1}{p}+\frac{1}{q}=1$. Then each 
solution $u(t,x)$ in Theorem A with $k=n+1$ has 
the following properties for all $\delta>0${\rm :} 
\[
\Vert \partial^{j}_t \partial^{\alpha}_x u(t,\cdot) \Vert_{L^q} \le 
C (1+\vert t \vert)^{-\left(\frac{n-1}{2}-\delta\right)
\left(\frac{1}{p}-\frac{1}{q} \right)} 
\sum_{i=0,1} \left\Vert u_i \right\Vert_{H^{N_p+j+|\alpha|-i,p}}
\]
where $N_p=\frac{3n+1}{2}\left(\frac{1}{p}-\frac{1}{q} \right)$, 
$j=0,1,2$, 
and $\alpha$ is any multi-index. 
\end{thm}

Based on Theorem \ref{thm:M}, we can develop the nonlinear scattering 
problems for the Kirchhof equation. But here, we want to 
exhibit the opposite 
phenomenon; for this, we will find the asymptotic profiles for the 
solutions to \eqref{Kirchhoff1}--\eqref{Kirchhoff2}. Let us present 
the definitions of free and non-free waves. 

\smallskip  

\noindent 
{\bf Definition.} {\em {\rm (i)} We say that $v_\pm=v_\pm(t,x)=\{v_+(t,x),
v_-(t,x)\}$ is a {\bf free wave} if it satisfies the equation 
$$\left(\partial^{2}_{t}-c^{2}_{\pm\infty}\Delta \right)v_\pm=0 
\quad \text{on $\mathbb{R} \times \mathbb{R}^n$.}$$ 

\noindent 
{\rm (ii)} Let $\sigma \ge 1$. We say that $v=v(t,x)$ is {\bf asymptotically 
free} in $\dot{H}^\sigma \times \dot{H}^{\sigma-1}$ if it is asymptotically
convergent to 
some free wave $v_\pm$ in $\dot{H}^\sigma \times \dot{H}^{\sigma-1}$, i.e., 
\[
\Vert v(t,\cdot)-v_\pm(t,\cdot)\Vert_{\dot{H}^\sigma}+
\Vert \partial_t v(t,\cdot)-\partial_t v_\pm(t,\cdot)
\Vert_{\dot{H}^{\sigma-1}}\to 0 \quad (t \to \pm\infty).
\] 

\noindent 
{\rm (iii)} Let $\sigma \ge 1$. We say that $w=w(t,x)$ is a {\bf 
non-free wave} in $\dot{H}^\sigma \times \dot{H}^{\sigma-1}$ if it is not 
asymptotically free.} \\

Theorem A states that each solution $u$ of 
\eqref{Kirchhoff1}--\eqref{Kirchhoff2} with initial data 
satisfying \eqref{smallness} with $k>2$, is asymptotically free. On the other 
hand, the next theorem states 
that the bound $k>2$ is sharp. More precisely, we have the following:

\begin{thm}[\cite{Matsuyama2}]\label{thm:thm2}
Assume that 
\[\text{either $n \ge 2$ and $1<k \le2$, or $n=1$ and $1<k<2$.}
\] 
Then there exists a solution $u(t,x)\in \cap_{s \ge 1} \dot{X}^s(\mathbb{R})$ 
of \eqref{Kirchhoff1}--\eqref{Kirchhoff2} with data satisfying 
\eqref{smallness}, which is a non-free wave in 
$\dot{H}^{\sigma}\times \dot{H}^{\sigma-1}$ for all $\sigma \ge1$. 
\end{thm}

The proof of Theorems \ref{thm:M}--\ref{thm:thm2} relies on the 
representation formulae for the corresponding linear equation. 
In \S3 we will introduce the representation formulae for more 
general strictly hyperbolic equations. 
Moreover, the argument of Theorem \ref{thm:thm2} 
is relating with the asymptotic 
behaviour of Bessel functions (see \cite{arons}). 


\section{Representation of solutions to linear Cauchy problems} 
In this section we introduce the representation formulae for more 
general equations than previously considered by using the 
asymptotic integration method along the argument of \cite{MR}. 
Let us consider the Cauchy problem for an $m^{\rm th}$ 
order strictly hyperbolic 
equation with time-dependent coefficients, for 
function $u=u(t,x)$: 
\begin{equation}
L(t,D_t,D_x)u 
\equiv D^{m}_t u+\sum_{\underset{j \le m-1}
{\vert \nu \vert+j=m}} 
a_{\nu,j}(t) D^{\nu}_x D^{j}_t u=0, \quad t \ne 0,
\label{Equation}
\end{equation}
with the initial condition 
\begin{equation}
D^{k}_t u(0,x)=f_k(x) \in C^{\infty}_0(\mathbb{R}^n), \quad 
k=0,1,\cdots,m-1, \quad x \in \mathbb{R}^n. 
\label{Initial condition}
\end{equation}
Denoting by $\mathcal{B}^{m-1}(\mathbb{R})$ the space of all 
functions whose derivatives up to $(m-1)^{\rm th}$ 
order are all bounded 
and continuous on $\mathbb{R}$, we assume that each 
$a_{\nu,j}(t)$ belongs to $\mathcal{B}^{m-1}(\mathbb{R})$ 
and satisfies 
\begin{equation}
\partial^{k}_t a_{\nu,j}(t) \in L^1(\mathbb{R}) \quad 
\text{for all $\nu,j$ with $|\nu|+j=m$, 
and $k=1,\ldots,m-1$}.
\label{L1} 
\end{equation}
Moreover, following the standard definition of equations
of the regularly hyperbolic type
(e.g. Mizohata \cite{Mizohata}), we will assume
that the symbol of the differential operator 
$L(t,D_t,D_x)$ has real and distinct roots 
$\varphi_1(t;\xi),\ldots,\varphi_m(t;\xi)$ 
for $\xi \ne0$, and 
\begin{equation}
L(t,\tau,\xi)=(\tau-\varphi_1(t;\xi)) \cdots 
(\tau-\varphi_m(t;\xi)), 
\label{strict hyperbolicity1}
\end{equation} 
\begin{equation}
\inf_{\underset{j \ne k}
{\vert \xi \vert=1,t \in \mathbb{R}}} 
\vert \varphi_j(t;\xi)-\varphi_k(t;\xi) \vert>0.
\label{strict hyperbolicity2}
\end{equation}

By applying the Fourier transform on $\mathbb{R}^n_x$ 
to \eqref{Equation}, 
we get 
\begin{equation}
D^{m}_t v+\sum^{m}_{j=1} h_j(t;\xi)D^{m-j}_{t} v=0, 
\label{wintner1}
\end{equation} 
where 
\[
h_j(t;\xi)=\sum_{\vert \nu \vert=j} a_{\nu,m-j}(t) \xi^\nu, 
\quad \xi \in \mathbb{R}^n.
\]
This is the ordinary differential equation, homogeneous 
of $m^{\rm th}$ order,
with the parameter $\xi=(\xi_1,\ldots,\xi_n)$. As usual, the strict 
hyperbolicity means that the characteristic roots of \eqref{wintner1} 
are real and can be written as 
$\varphi_1(t;\xi),\ldots,\varphi_m(t;\xi)$ satisfying 
\eqref{strict hyperbolicity1}--\eqref{strict hyperbolicity2}. 
Notice that each 
$\varphi_\ell(t;\xi)$ has a homogeneous degree one with respect to 
$\xi$. In this section we will establish the representation formulae
for solutions of the Cauchy problem \eqref{Equation} in the
form of the oscillatory integrals. Let $\widehat{u}(t;\xi)$ be the solution of 
\eqref{wintner1} with the intial data $\widehat{f}_k(\xi)$ 
$(k=0,\ldots,m-1)$. Let $v_k(t;\xi)$ be the solution of \eqref{wintner1} 
with $(D^{j}_t v_k)(0;\xi)=\delta^{j}_k$ for $j,k=0,1,\ldots,m-1$. We 
set 
\[W(t;\xi)=\left(
\begin{array}{cccc}
v_0(t;\xi) & v_1(t;\xi) & \cdots & v_{m-1}(t;\xi) \\
D_t v_0(t;\xi) & D_t v_1(t;\xi) & \cdots & D_t v_{m-1}(t;\xi) \\
\vdots & \vdots & \cdots & \vdots \\ 
D^{m-1}_t v_0(t;\xi) & D^{m-1}_t v_1(t;\xi) & \cdots & D^{m-1}_t v_{m-1}(t;\xi)
\end{array}\right). 
\]
Hence, $W(t;\xi)$ is the fundamental matrix of \eqref{wintner1}. 
Defining 
$$\vartheta_j(t;\xi)=\int\limits^{t}_{0} \varphi_j(s;\xi) \, ds, 
\quad j=1,\ldots,m,$$ 
we introduce the matrix 
\[Y(t;\xi)=
\begin{pmatrix} 
\mathrm{e}^{i\vartheta_1(t;\xi)} & \cdots 
& \mathrm{e}^{i\vartheta_m(t;\xi)}\\ 
D_t \mathrm{e}^{i\vartheta_1(t;\xi)} & \cdots & D_t 
\mathrm{e}^{i\vartheta_m(t;\xi)}\\ 
\cdots & \cdots & \cdots\\ 
D^{m-1}_t \mathrm{e}^{i\vartheta_1(t;\xi)} & \cdots & D^{m-1}_t 
\mathrm{e}^{i\vartheta_m(t;\xi)}
\end{pmatrix}. 
\]
Matrix $Y(t;\xi)$ is the fundamental matrix of 
a perturbed ordinary differential 
equation of \eqref{wintner1}: 
\[
\left(D_t-\varphi_1(t;\xi) \right)\cdots \left(D_t-\varphi_m(t;\xi) 
\right)w=0. 
\]
Then we can write this equation as 
\begin{equation}
D^{m}_t w+\sum^{m}_{j=1}h_j(t;\xi)D^{m-j}_{t}w
+\sum^{m}_{j=2}\widetilde{h}_j(t;\xi)D^{m-j}_{t}w=0, 
\label{perturbation}
\end{equation}
where $\widetilde{h}_j(t;\xi)$ satisfies 
\[
\widetilde{h}_j(t;\xi)=\begin{cases} 
\quad 0, & \quad j=1,\\ 
\displaystyle{\sum_{\underset{(\nu_2,\ldots,\nu_j)\ne(0,\ldots,0)}{1 \le \vert 
\nu \vert \le j-1}}}
\widetilde{c}_{\nu_1 \ldots\nu_j}\;\varphi^{\nu_1}_{\ell_1} \;
(D_t \varphi_{\ell_2})^{\nu_2} \cdots (D^{m-j+1}_t \varphi_{\ell_j})^{\nu_j}, 
& \quad j=2,\ldots,m, 
\end{cases}
\]
with some constants $\widetilde{c}_{\nu_1 \ldots\nu_j} \ne 0$. This means that 
each $\mathrm{e}^{i\vartheta_\ell(t;\xi)}$ satisfies \eqref{perturbation}, and 
$\mathrm{e}^{i\vartheta_1(t;\xi)},\ldots,\mathrm{e}^{i\vartheta_m(t;\xi)}$ 
are linearly independent for $\xi \ne 0$ and $t \in \mathbb{R}$. It can be 
checked that the coefficient $\widetilde{h}_1(t;\xi)$ of $D^{m-1}_t w$ always 
vanishes for every $m$, by an induction argument on $m$. 
Then it 
follows from Proposition 2.4 of \cite{MR} (cf. \cite{Ascoli,Wintner}) 
that there exists the limit
\begin{equation}
\displaystyle{\lim_{t \to \pm\infty}}Y(t;\xi)^{-1}W(t;\xi)=L_\pm(\xi). 
\label{a1}
\end{equation}
Set 
\[
L_\pm(\xi)=\left(
\begin{array}{cccc}
\alpha^{1}_{0,\pm}(\xi) & \alpha^{1}_{1,\pm}(\xi) & \cdots & 
\alpha^{1}_{m-1,\pm}(\xi) \\
\alpha^{2}_{0,\pm}(\xi) & \alpha^{2}_{1,\pm}(\xi) & \cdots & 
\alpha^{2}_{m-1,\pm}(\xi) \\
\vdots & \vdots &\cdots &\vdots \\
\alpha^{m}_{0,\pm}(\xi) & \alpha^{m}_{1,\pm}(\xi) & \cdots & 
\alpha^{m}_{m-1,\pm}(\xi) 
\end{array}\right).
\]
Furthermore, writing 
\begin{eqnarray}
&& R_\pm(t;\xi)=Y(t;\xi)^{-1}W(t;\xi)-L_\pm(\xi) \label{e1}\\
&=& \left(
\begin{array}{cccc}
\varepsilon^{1}_{0,\pm}(t;\xi) & \varepsilon^{1}_{1,\pm}(t;\xi) & \cdots & 
\varepsilon^{1}_{m-1,\pm}(t;\xi) \\
\varepsilon^{2}_{0,\pm}(t;\xi) & \varepsilon^{2}_{1,\pm}(t;\xi) & \cdots & 
\varepsilon^{2}_{m-1,\pm}(t;\xi) \\
\vdots & \vdots &\cdots &\vdots \\
\varepsilon^{m}_{0,\pm}(t;\xi) & \varepsilon^{m}_{1,\pm}(t;\xi) & \cdots & 
\varepsilon^{m}_{m-1,\pm}(t;\xi) 
\end{array}\right), \nonumber 
\end{eqnarray}
we have 
\[
W(t;\xi)=Y(t;\xi)\left(L_\pm(\xi)+R_\pm(t;\xi)\right). 
\]
Thus we arrive at 
\begin{eqnarray*}
D^{\ell}_t v_k(t;\xi)&=& \sum^{m}_{j=1} \left(\alpha^{j}_{k,\pm}(\xi)
+\varepsilon^{j}_{k,\pm}(t;\xi)\right) D^{\ell}_t 
\mathrm{e}^{i\vartheta_j(t;\xi)}\nonumber \\
&=& \sum^{m}_{j=1} \left(\alpha^{j}_{k,\pm}(\xi)
+\varepsilon^{j}_{k,\pm}(t;\xi)\right) p_\ell(\varphi_j(t;\xi))
\mathrm{e}^{i\vartheta_j(t;\xi)}
\end{eqnarray*}
for $k,\ell=0,\cdots,m-1$, 
where each $p_\ell(\varphi_j(t;\xi))$ 
is determined by the equation $$D^{\ell}_t \mathrm{e}^{i\vartheta_j(t;\xi)}
=p_\ell(\varphi_j(t;\xi))\mathrm{e}^{i\vartheta_j(t;\xi)}.$$

We note that for the second order equations we have $m=2$ and 
the next theorem covers the case of the wave equation as a 
special case, 
also improving the corresponding result in 
\cite{Matsuyama1,Matsuyama0,Matsuyama2}. 
The result is as follows: 
\begin{thm} \label{thm:asymptotic}
Assume that the characteristic roots 
$\varphi_1(t;\xi),\ldots,\varphi_m(t;\xi)$ of \eqref{wintner1} are 
real and distinct for all $t \in \mathbb{R}$ and for all 
$\xi \in \mathbb{R}^n\backslash0$, and that they
satisfy \eqref{strict hyperbolicity2}.
Then there exists 
$\alpha^{j}_{k,\pm}(\xi)$ and $\varepsilon^{j}_{k,\pm}(t;\xi)$ 
determined by \eqref{a1} and \eqref{e1}, 
respectively, such that the solution 
$u(t,x)$ of our problem 
$\eqref{Equation}$--$\eqref{Initial condition}$ is 
represented by 
\[
D^{\ell}_t u(t,x)=\sum^{m-1}_{k=0}\sum^{m}_{j=1} 
\mathcal{F}^{-1} 
\left[ \left(\alpha^{j}_{k,\pm}(\xi)
+\varepsilon^{j}_{k,\pm}(t;\xi)\right) p_\ell(\varphi_j(t;\xi)) 
\mathrm{e}^{i\vartheta_j(t;\xi)} \widehat{f}_k (\xi)\right](x), 
\quad t \gtrless 0,
\] 
for $\ell =0,\ldots,m-1$, where 
\[
\left\vert \alpha^{j}_{k,\pm}(\xi)\right\vert 
\le c\vert \xi \vert^{-k}, 
\quad \left\vert \varepsilon^{j}_{k,\pm}(t;\xi)\right\vert 
\le c\vert \xi \vert^{-k} 
\int\limits^{+\infty}_{\vert t \vert} \Psi(s) \, ds, 
\] 
and $\Psi(t)$ is given by 
\[
\Psi(t)=\sum_{\underset{j\le m-2}{\vert \nu \vert+j=m}} 
\left\vert \partial_t a_{\nu,j}(t)\right\vert \cdots 
\left\vert \partial^{m-j-1}_t a_{\nu,j}(t)\right\vert.
\]
For the higher order derivatives of 
amplitude functions, we have, for $|\mu|\ge 1$, 
\[
\left\vert D^{\mu}_{\xi}\alpha^{j}_{k,\pm}(\xi)\right\vert 
\le c |\xi|^{-k}, \quad 
\left\vert D^{\mu}_{\xi} \varepsilon^{j}_{k,\pm}(t;\xi)\right\vert  
\le c \mathrm{e}^{\int\limits^{|t|}_{0} 
(1+s)^{\vert \mu \vert}\Psi(s) \, ds}|\xi|^{-k}, \quad 
\quad |\xi| \ge 1,
\]
\[
\left\vert D^{\mu}_{\xi}\alpha^{j}_{k,\pm}(\xi)\right\vert 
\le c |\xi|^{-k-|\mu|}, \quad 
\left\vert D^{\mu}_{\xi} \varepsilon^{j}_{k,\pm}(t;\xi)\right\vert  
\le c \mathrm{e}^{\int\limits^{|t|}_{0} 
(1+s)^{\vert \mu \vert}\Psi(s) \, ds}|\xi|^{-k-|\mu|}, \quad 
0<|\xi|<1.
\]

If we further assume 
that 
$$(1+\vert t \vert)^{\vert \mu \vert}\partial^{k}_t a_{\nu,j}
(t;\xi) \in 
L^1(\mathbb{R})$$ for some $\mu$ with $\vert \mu \vert \ge 1$, 
and for all 
$\nu,j$, and $k=1,\ldots,m-1$, then the bound of each 
$D^{\mu}_{\xi} \varepsilon^{j}_{k,\pm}(t;\xi)$ is uniform in $t$.
\end{thm}

\end{document}